\documentclass[11pt]{article}
\usepackage{amsmath,amsthm,amssymb}
\newcommand{\remove}[1]{}
\newtheorem{ithm}{Theorem}
\newtheorem{thm}{Theorem}[section]
\newtheorem{claim}[thm]{Claim}
\newtheorem{lem}[thm]{Lemma}
\newtheorem{cor}[thm]{Corollary}

\newtheorem{fact}[thm]{Fact}
\newtheorem{prop}[thm]{Proposition}

\theoremstyle{definition}
\newtheorem{define}[thm]{Definition}

\newtheorem{remark}[thm]{Remark}

\def\F{{\mathbb{F}}}
\def\R{{\mathbb{R}}}
\def\Z{{\mathbb{Z}}}

\def\C{{\mathbb{C}}}

\def\L{{\mathcal L}}

\def\SC{{\mathbf{SC}}}

\def\fl{\text{fl}}
\def\SG{\text{SG}}

\newcommand{\ip}[2]{\langle #1,#2 \rangle}

\def\_{\,\,\,\,\,}

\def\rank{\text{rank}}
\def\span{\text{span}}
\def\dimension{\text{dim}}

\def\poly{\text{poly}}
\def\tr{\text{tr}}

\def\poly{\text{poly}}

\newcommand{\eps}{\epsilon}

\DeclareMathOperator{\E}{\mathbb{E}}

\begin{document}

\title{Rank Bounds for Design Matrices  with Applications  to
Combinatorial Geometry and Locally Correctable Codes}
\author{Boaz Barak\thanks{Microsoft Research New England.
Email: \texttt{boaz@microsoft.com}. Most of the work done while at Princeton University and supported by
NSF grants CNS-0627526, CCF-0426582 and CCF-0832797, and the Packard
and Sloan fellowships.} \and
Zeev Dvir\thanks{Department of Computer Science, Princeton University.
Email: \texttt{zeev.dvir@gmail.com}. Research partially supported by NSF grant CCF-0832797 and by the Packard fellowship.} \and
Avi Wigderson\thanks{School of Mathematics, Institute for Advanced Study.
Email: \texttt{avi@ias.edu}. Research partially supported by NSF grants CCF-0832797 and DMS-0835373.} \and
Amir Yehudayoff\thanks{Department of Mathematics,
Technion - IIT.
Email: \texttt{amir.yehudayoff@gmail.com}. Most of the work done while at the IAS. Research partially supported by NSF grants CCF-0832797 and DMS-0835373.}}
\date{}
\maketitle

\begin{abstract}
A $(q,k,t)$-design matrix is an $m \times n$ matrix whose pattern of zeros/non-zeros satisfies the following  design-like
condition: each row has at most $q$ non-zeros, each column has at least $k$ non-zeros and the supports of every two columns
intersect in at most $t$ rows. We prove that for $m\geq n$, the rank of any $(q,k,t)$-design matrix over a field of
characteristic zero (or sufficiently large finite characteristic) is at least \[ n - \left(\frac{qtn}{2k}\right)^2 . \]
Using this result we derive the following applications:
\begin{description}
\item[Impossibility results for $2$-query LCCs over large fields.]
A $2$-query locally correctable code (LCC) is an error correcting code  in which every codeword coordinate can be
    recovered, probabilistically, by reading at most two other code positions.  Such codes have numerous applications and
    constructions (with exponential encoding length) are known over finite fields of small characteristic. We show that
    infinite families of such linear $2$-query LCCs \emph{do not exist} over fields of
    characteristic zero or large characteristic regardless of the encoding length.

\item[Generalization of known results in combinatorial geometry.]
We prove a quantitative analog of the Sylvester-Gallai theorem:
Let $v_1,\ldots,v_m$ be a set of points in $\C^d$ such
that for every $i \in [m]$ there exists at least $\delta m$ values of $j \in [m]$ such that the line through $v_i,v_j$
contains a third point in the set. We show that the dimension of $\{ v_1,\ldots,v_m \}$ is at most $O(1/\delta^2)$. Our results generalize to the high dimensional case (replacing lines with planes, etc.) and to the case where the points are colored (as in the Motzkin-Rabin Theorem).
\end{description}
\end{abstract}

\newpage
\section{Introduction}\label{sec-intro}

In this work we study what \emph{combinatorial} properties of matrices guarantee high algebraic \emph{rank}, where a property is
combinatorial if it depends only on the zero/non-zero pattern of the matrix, and not on the values of its entries.  This question has a rich history in mathematics (see Section~\ref{sec: rel work}), and some computer science motivations: 

\begin{description}

\item[Locally correctable codes.] A \emph{locally correctable code} is an error correcting code in which for every
    codeword $y$, given a corrupted version $\Tilde{y}$ of $y$ and an index $i$, one can recover the correct value of $y_i$ from
    $\Tilde{y}$ by looking only at very few coordinates of $\Tilde{y}$. It is an open question in coding theory to understand
    the tradeoffs between the fraction of errors, locality (number of coordinates read) and rate (ratio of message length
    to codeword length) of such codes, with very large gaps between the known upper bounds and lower bounds (see the survey
    \cite{Tre04}). The question is open even for linear codes, where the condition of being locally correctable turns out to
    be equivalent to the existence of low weight codewords in the dual codewords that are ``well-spread'' in some precise
    technical sense (see Section~\ref{sec-LCC}). Because of the relation between the rate of the code and its dual, the
    question becomes equivalent to asking whether this combinatorial ``well-spreadness'' condition guarantees high rank.

\item[Matrix rigidity.] A longstanding question is to come up with an explicit matrix that is \emph{rigid} in the sense
    that its rank cannot be reduced by changing a small number of its entries. Random matrices are extremely rigid, and
    sufficiently good explicit constructions will yield lower bounds for arithmetic circuits~\cite{Val77}, though we are
    still very far from achieving this (see the survey~\cite{Lok09}). One can hope that a combinatorial property guaranteeing
    large rank will be robust under small perturbations, and hence a matrix satisfying such a property will automatically be
    rigid.

\end{description}

In both these cases it is crucial to obtain bounds on the rank that depend solely on the zero/non-zero pattern of the matrix,
without placing any restrictions on the non-zero coefficients. For example, there are very strong bounds known for matrix
rigidity under the restriction that the non-zero coefficients have bounded magnitude (see Chapter 3 in \cite{Lok09}), but they only imply lower bounds
in a very restricted model. In fact, there is a relation between the two questions, and sufficiently good answers for the first
question will imply answers for the second one~\cite{Dvi10}. We stress that these two examples are in no way exhaustive. The
interplay between combinatorial and algebraic properties of matrices is a fascinating question with many potential applications
that is still very poorly understood.

\subsection{Our Results}

%
%

In this work we give a combinatorial property of complex matrices that implies high rank. While not strong enough to prove
rigidity results, we are able to use it to obtain several applications in combinatorial geometry and locally correctable codes.
Our main result is the following theorem, giving a lower bound on the rank of matrix whose non-zero pattern forms has  certain combinatorial-design like properties in the sense that the sets of non-zero entries in each column have small intersections. (This theorem is restated as Theorem~\ref{thm-rankdesign}.)

\begin{ithm}[Rank bound for design matrices] \label{ithm:main}
Let $m\geq n$. We say that an  $m\times n$ complex matrix $A$ is a $(q,k,t)$-design matrix if every row of $A$ has at most $q$
non-zero entries, every column of $A$ has at least $k$ non-zeroes entries, and the supports of every two columns intersect in at most $t$ rows. For every such $A$,
\[ \rank(A) \geq n - \left(\frac{q\cdot t \cdot n}{2k}\right)^2 . \]
\end{ithm}

We also show that Theorem~\ref{ithm:main}, and in fact any result connecting the zero/non-zero pattern to rank, can be made to
hold over arbitrary characteristic zero fields and also over fields of sufficiently large (depending on $m,n$) finite
characteristic.

\subsubsection{Applications to Combinatorial Geometry}

Our most immediate applications of Theorem~\ref{ithm:main} are to questions regarding line-point incidences. Results on
line-point incidences have recently found use in the area of computational complexity in relation to pseudo-randomness
\cite{BKT04,BIW06} and de-randomization \cite{KS09, SS10}. In this setting we have an arrangement of a finite number of points in
real or complex space. Every  such arrangement  gives rise to a set of lines, namely, those lines that pass through at least two
of the points in the arrangement. Information about these lines can be converted, in some cases, into information about the
dimension of the set of points (i.e. the dimension of the space the points span). Our rank theorem can be used to derive
generalizations for two well-known theorems in this area: the Sylvester-Gallai theorem and the Motzkin-Rabin theorem.

\paragraph{Generalizing the Sylvester-Gallai Theorem.}
The Sylvester-Gallai (SG for short) theorem says that if $m$ distinct points $v_1,\ldots,v_m \in \R^d$ are not collinear,
then there exists a line that passes through exactly two of them. In its contrapositive form the SG theorem says that if for every $i
\neq j$ the line through $v_i$ and $v_j$ passes through a third point $v_k$, then $\dim\{v_1,\ldots,v_m \} \leq 1$,
where $\dim\{v_1,\ldots,v_m \}$ is the dimension of the smallest affine subspace containing the points. This theorem was first conjectured by Sylvester in 1893 \cite{Sylvester}, proved (in dual form) by Melchior in 1940 \cite{Melchior},
and then independently conjectured by Erdos in 1943 \cite{Erdos43} and proved by Gallai in 1944. The SG theorem has several beautiful proofs and
many generalizations, see the survey \cite{BM90}. Over the complex numbers the (tight) bound on the dimension is  $2$ instead of $1$. The
complex version was first proven by Kelly \cite{Kel86} using a deep results from algebraic geometry, and more recently, an
elementary proof was found by Elkies, Pretorius and Swanepoel \cite{EPS06} who also proved it over the quaternions with an upper of 4 on the dimension.

We say that the points $v_1,\ldots,v_m$ (in $\R^d$ or $\C^d$) form a $\delta$-SG configuration if for every $i \in [m]$ there exists at least $\delta m$ values of $j \in [m]$ such that the line through $v_i,v_j$ contains a third point in the set. Szemeredi and Trotter~\cite{ST83} showed that, when $\delta$ is larger than some absolute constant close to 1, then the dimension of a $\delta$-SG configuration is at most one (over the reals). We show the following generalization of
their result to arbitrary $\delta>0$ (and over the complex numbers).

\begin{ithm}[Quantitative SG theorem]\label{ithm:SG}
If $v_1,\ldots,v_m \in \C^d$ is a $\delta$-SG configuration then $\dim\{ v_1,\ldots,v_m \} < 13/\delta^2$.
\end{ithm}

We note that one cannot replace the bound $13/\delta^2$ of Theorem~\ref{ithm:SG} with $1$ or even with any fixed constant, as one
can easily create a $\delta$-SG configuration of dimension roughly $2/\delta$ by placing the points on $1/\delta$ lines. This is
analogous to error correcting codes, where once the fraction $\delta$ of agreement between the original and corrupted codeword
drops below half there can be no unique decoding. In that sense our result can be thought of as a \emph{list decoding} variant of
the SG theorem, whereas the result of \cite{ST83} is its unique decoding variant. We also show an ``average case'' version of the SG
theorem, proving a bound on the dimension of a large subset of the points under the assumption that there are many collinear triples
(see Theorem~\ref{thm-SGav}).

We also prove a version of Theorem~\ref{thm-deltaSG} with lines replaced by $k$-flats ($k$-dimensional affine subspaces). This generalizes a theorem of Hansen \cite{Han65,BE67} which deals with the case $\alpha=1$. The statement of this result is technical
and so we give it in Section~\ref{sec-highdim} where it is also proven.

Since our proofs use elementary (and purely algebraic) reductions to the rank theorem, they hold over arbitrary fields of
characteristic zero or of sufficiently large finite characteristic. This is in contrast to many of the known proofs of such
theorems which often rely on specific properties of the real (or complex) numbers. However, we currently do not recover the full
version of the original SG theorem, in the sense that even for $\delta=1$ we do not get a bound of $1$ (or $2$ for complex
numbers) on the dimension. (However, the term $13/\delta^2$ can be improved a bit in the $\delta=1$ case to obtain a bound of
$9$ on the dimension.)

\paragraph{Generalizing the Motzkin-Rabin Theorem.} The Motzkin-Rabin (MR for short) theorem (see e.g. \cite{BM90})
is an interesting variant of the Sylvester-Gallai theorem that states that if points $v_1,\ldots,v_m \in \R^d$ are colored
either red or blue and there is no monochromatic line passing through at least two points, then they are all collinear. As in the SG
theorem, we obtain a quantitative generalization of the MR theorem such that (letting $b$ and $r$ be the numbers of blue and red
points respectively), if for every blue (resp. red) point $v$,  there are $\delta b$ blue (resp. $\delta r$ red) points $v'$
where the line through $v$ and $v'$ passes through a red (resp. blue) point, then $\dim\{ v_1,\ldots,v_m \} \leq
O(1/\delta^4)$. We also prove a three colors variant of the MR theorem, showing that if $v_1,\ldots,v_m$ are colored red, blue and
green, and all lines are not monochromatic, then $\dim\{v_1,\ldots,v_m\}$ is at most some absolute constant.

\subsubsection{Locally Correctable Codes}

A (linear) $q$ query locally correctable code ($(q,\delta)$-LCC for short) over a field $\F$ is a subspace $C \subseteq \F^n$
such that, given an element $\Tilde{y}$ that disagrees with some $y \in C$ in at most $\delta n$ positions and an index $i\in
[n]$, one can recover $y_i$  with, say, probability $0.9$, by reading at most $q$ coordinates of $\Tilde{y}$. Over the field of
two elements $\F_2$ the standard Hadamard code construction yields a (2,$\delta$)-query LCC with dimension $\Omega(\log(n))$ for
constant $\delta>0$ (see the survey \cite{Tre04}). In contrast we show that for every constant $\delta>0$ there do not exist
infinite family of such codes over the complex numbers:

\begin{ithm}[Impossibility of $2$-query LCCs over $\C$] \label{ithm:lcc}
If $C$ is a $2$-query LCC for $\delta$ fraction of errors over $\C$, then $\dim( C ) \leq O(1/\delta^9)$.
\end{ithm}

We note that the Hadamard construction does yield a \emph{locally decodable code} over
the complex numbers with dimension $\Omega(\log n)$. Locally decodable codes are the  relaxation of a locally correctable codes
where one only needs to be able to recover the coordinates of the original message as opposed to the codeword. Thus over the
complex numbers, there is a very strong separation between the notions of locally decodable and locally correctable codes,
whereas it is consistent with our knowledge that for, say, $\F_2$ the rate/locality tradeoffs  of both notions are the same.

\subsection{Related Work}
\label{sec: rel work}
The idea to use matrix scaling  to study structural properties of matrices was already present in \cite{CPR00}. This work, which was also motivated by the problem of matrix rigidity, studies the presence of short cycles in the graphs of non-zero entries of a square matrix.

A related line of work on the rank of `design' matrices is the work emerging from Hamada's conjecture \cite{Ham73}. (See
\cite{JT09} for a recent result and more references.) Here, a design matrix is defined using stricter conditions (each row/column
has exactly the same number of non-zeros and the intersections are also all of the same size) which are more common in the
literature dealing with combinatorial designs. In order to be completely consistent with this line of work we should have called
our matrices `approximate-design' matrices. We chose to use the (already overused) word `design' to make the presentation more
readable. We also note that considering approximate designs only makes our results stronger. Hamada's conjecture states that of
all zero/one matrices whose support comes from a design (in the stricter sense), the minimal rank is obtained by matrices coming
from geometric designs (in our language, Reed-Muller codes). In contrast to this paper, the emphasis in this line of works is typically on small finite fields.  We
note here that the connection between Hamada's conjecture and LCCs was already observed by Barkol, Ishai and Weinreb \cite{BIW07}
who also conjectured (over small fields) the `approximate-design' versions which we prove here for large fields.

Another place where the support of a matrix is connected to its rank is in graph theory where we are interested in minimizing the
rank of a (square, symmetric) real matrix which has the same support as the adjacency matrix of a given graph. This line of work
goes back for over fifty years and has many applications in graph theory. See \cite{FH07} for a recent survey on this topic.

Over the reals we can also ask about the minimal rank of matrices with certain {\em sign-pattern}. That is, given a matrix over
$\{1,-1\}$, what is the minimal rank of a matrix which has the same sign-pattern. This minimal rank is called the {\em sign-rank}
of a matrix. The question of coming up  with (combinatorial or otherwise) properties that imply high sign-rank is one of major
importance and has strong connections to  communication complexity, learning theory and circuit complexity, among others. For a
recent work with plenty of references see \cite{RS08}. In particular we would like to mention a connection to the work of Forster
\cite{For02} on the sign-rank of the Hadamard matrix. (An earlier version of this work used a variant \cite{Bar98,Har10} of
a lemma from \cite{For02} instead of the results of \cite{RS89} on matrix scaling to obtain our main result.)

\subsection{Organization}

In Section~\ref{sec-techniques} we give a high level overview of our techniques. In Section~\ref{sec-rank} we prove our main
result on the rank of design matrices. In Section~\ref{sec-SG} we prove our quantitative variants of the Sylvester-Gallai
theorem. In Section~\ref{sec-highdim} we prove the high-dimensional analog of Theorem~\ref{thm-deltaSG} where lines are replaced with flats. In Section~\ref{sec-MR} we prove our generalizations of the Motzkin-Rabin theorem. In Section~\ref{sec-LCC} we prove our
results on locally correctable codes. In Section~\ref{sec-finite}  we show how our results extend to other fields. We conclude in
Section~\ref{sec-open} with a discussion of open problems.

\section{Our Techniques} \label{sec-techniques}

We now give high-level proof overviews for some of our results.

\subsection{Rank Lower Bounds for Design Matrices}

Theorem~\ref{ithm:main} -- the rank lower bound for design matrices -- is proved in two steps.  We now sketch the proof,
ignoring some subtleties and optimizations. The proof starts with the observation that, as in the case of matrix rigidity and
similar questions, the result is much easier to prove given a bound on the \emph{magnitude} of the non-zero entries. Indeed, if
$A$ is a $(q,k,t)$-design matrix and all of its non-zero entries have absolute value in $[1/c,1]$ for some constant $c$, then the $n\times n$ matrix $M=A^*A$ is  \emph{diagonally dominant}, in the sense that for all $i\neq j$, $m_{ii} \geq k/c^2$ but $|m_{ij}| \leq t$. (Here $A^*$ denotes the conjugate transpose of $A$.) Thus one can use known results on such
matrices (e.g. \cite{Alo09}) to argue that $\rank(A) \geq \rank(M) \geq n - (ntc^2/k)^2$. Our main idea is to reduce to this case
where the non-zero coefficients of $A$ are (roughly) bounded using \emph{matrix scaling}.

A \emph{scaling} $\Hat{A}$ of a matrix $A$ is obtained by multiplying for all $i,j$,
the $i$'th row of $A$ by some positive number $\rho_i$ and
the $j$'th column of $A$ by some positive number $\gamma_j$.
Clearly, $A$ and $\Hat{A}$ share the same rank and zero/non-zero pattern.
We use known matrix-scaling results \cite{Sinkhorn,RS89}
to show that every $(q,k,t)$-design
matrix $A$ has a scaling in which every entry has magnitude at most (roughly) $1$ but its columns have norm at least (roughly) $\sqrt{k/q}$.
We note that the typical application of matrix-scaling was
with respect to the $\ell_1$-norm of the rows and columns.
Here we take a different path: We use scaling with respect to $\ell_2$-norm.

We defer the description of this step to Section~\ref{sec-rank} but the high level idea is to use a theorem of \cite{RS89} that
shows that such a scaling exists (in fact without the dependence on $q$) if $A$ had the property of not containing any large all-zero sub-matrix. While this property cannot be in general guaranteed, we show that by repeating some rows of  $A$ one can obtain a
matrix $B$ that has this property, and a scaling of $B$ can be converted into a scaling of $A$.  Since our lower bound on the
entry $m_{ii}$ in the bounded coefficient case (where again $M=A^*A$) only used the fact that the columns have large norms, we
can use the same argument as above to lower bound the rank of $M$, and hence of $A$.

\subsection{Generalized Sylvester-Gallai Theorem}

Recall that the quantitative SG theorem (Theorem~\ref{ithm:SG}) states that every $\delta$-SG configuration $v_1,\ldots,v_n$, has
dimension at most $13/\delta^2$. Our proof of Theorem~\ref{ithm:SG} uses Theorem~\ref{ithm:main} as
follows. Suppose for starters that every one of these lines passed through \emph{exactly} three points. Each such line induces an
equation of the form $\alpha v_i + \beta v_j + \gamma v_k = 0$. Now for $m=\delta n^2$, let $A$ be the $m \times n$ matrix whose
rows correspond to these equations. Since every two points participate in only one line, $A$ will be a $(3,\delta n, 1)$ design matrix, meaning that according to Theorem~\ref{ithm:main}, $A$'s rank is at least $n-\left(\tfrac{3}{2\delta}\right)^2$.
Since $A$ times the matrix whose rows are $v_1,\ldots,v_n$ is zero we have
$\dim\{ v_1,\ldots,  v_n \} \leq n - \rank(A)$.
We thus get an upper bound of $\lfloor 9/4 \rfloor=2$ on this dimension. To handle
the case when some lines contain more than three points, we choose in some careful way from each line $\ell$ containing $r$
points a subset of the $\binom{r}{3}$ equations of the form above that it induces on its points. We show that at some small loss
in the parameters we can still ensure the set of equations forms a design, hence again deriving a lower bound on its rank via
Theorem~\ref{ithm:main}.

Our method extend also to an ``average case'' SG theorem (Theorem~\ref{thm-SGav}), where one only requires that the set of points
supports many (i.e., $\Omega(n^2)$) collinear triples and that each pair of points appear together in a few collinear triples. In this case we are able to show that there is a subset of $\Omega(n)$
points whose span has dimension $O(1)$. See Section~\ref{sec-SG} for more details. Our generalizations  of  the Motzkin-Rabin theorem follow from our theorem on $\delta$-SG configurations via simple reductions (see Section~\ref{sec-MR}).

\subsection{Locally Correctable Codes}
\label{ssec: lcc}

At first sight, Theorem~\ref{ithm:lcc} -- non existence of $2$ query locally correctable codes over $\C$ -- seems like it should
be an immediate corollary of  Theorem~\ref{ithm:SG}. Suppose that a code $C$ maps $\C^d$ to $\C^n$, and let $v_1,\ldots,v_n$
denote the rows of its generating matrix. That is, the code maps a message $x \in \C^d$ to the vector
$(\ip{v_1}{x},\ldots,\ip{v_n}{x})$. The fact that $C$ is a $2$ query LCC for $\delta$ errors implies that for every such row
$v_i$, there are roughly $\delta n$ pairs $j,k$ such that $v_i$ is in the span of $\{ v_j,v_k\}$. Using some simple scaling/change of basis, this gives precisely the condition of  being a
$\delta$-SG configuration, save for one caveat: In a code there is no guarantee that all the vectors $v_1,\ldots,v_n$ are
distinct. That is, the code may have repeated coordinates that are always identical. Intuitively it seems that such repetitions
should not help at all in constructing LCCs but proving this turned out to be elusive. In fact, our proof of
Theorem~\ref{ithm:lcc} is rather more complicated than the proof Theorem~\ref{ithm:SG}, involving repeated applications of
Theorem~\ref{ithm:main} which result also in somewhat poorer quantitative bounds. The idea behind the proof to use a variant of
the ``average case'' SG theorem to repeatedly find $\Omega(n)$ points among $v_1,\ldots,v_n$ whose span has $O(1)$ dimension,
until there are no more points left. We defer all details to Section~\ref{sec-LCC}.

Given Theorem~\ref{ithm:main}, one may have expected that Theorem~\ref{ithm:lcc} could be extended for LCCs of any constant number $q$ of queries. After all, the condition of $C$ being an LCC intuitively seems like only a slight relaxation of requiring
that the dual code of $C$ has a generating matrix whose non-zero pattern is a combinatorial design, and indeed in known constructions of LCCs, the dual code does form a design.
We are not, however, able to extend our results to $3$ and more queries.
A partial explanation to our inability is that $3$ query LCCs give rise to configuration of planes (instead of lines) and point
and planes exhibit much more complicated combinatorial properties than lines.

\section{Rank of Design Matrices}\label{sec-rank}

In this section we prove our main result which gives a lower bound on the rank of matrices whose zero/non-zero pattern satisfies certain properties. We start by defining these properties formally.

\begin{define}[Design matrix]\label{def-designmatrix}
Let $A$ be an $m \times n$ matrix over some field. For $i \in [m]$ let $R_i \subset [n]$ denote the set of indices of all non-zero entries in the $i$'th row of $A$. Similarly, let $C_j \subset [m]$, $j \in [n]$,
denote the set of non-zero indices in the $j$'th column. We say that $A$ is a {\em $(q,k,t)$-design matrix} if
\begin{enumerate}
\item For all $i \in [m]$, $|R_i| \leq q$.
\item For all $j \in [n]$, $|C_j| \geq k$.
\item For all $j_1 \neq j_2 \in [n]$, $|C_{j_1} \cap C_{j_2} | \leq t$.
\end{enumerate}
\end{define}

\begin{thm}[Restatement of Theorem~\ref{ithm:main} -- rank of design matrices]\label{thm-rankdesign}
Let $A$ be an $m \times n$ complex matrix. If $A$ is a $(q,k,t)$-design matrix then
$$\rank(A) \geq n - \left(\frac{q\cdot t \cdot n}{2k}\right)^2 .$$
\end{thm}

\begin{remark}
The proof of the theorem actually holds under a slightly weaker condition on the sizes of the intersections. Instead of requiring that $|C_{j_1} \cap C_{j_2} | \leq t$ for all pairs of columns $j_1 \neq j_2$, it is enough to ask that $$ \sum_{j_1 \neq j_2} |C_{j_1} \cap C_{j_2} |^2 \leq n^2 \cdot t^2. $$ That is, there could be some pairs with large intersection as long as the average of the squares is not too large.
\end{remark}

The proof of the theorem is given below, following some preliminaries.

\subsection{Preliminaries for the Proof of Theorem~\ref{thm-rankdesign}}

\paragraph{Notation:} For a set of real vectors $V \in \C^n$ we denote by $\rank(V)$ the dimension of the vector space spanned by elements of $V$. We denote the $\ell_2$-norm of a vector $v$ by $\| v \|$.
We denote by $I_n$ the  $n \times n$ identity matrix.

We start with definitions and results on matrix scaling.

\begin{define}\label{def-scaling}[Matrix scaling]
Let $A$ be an $m \times n$ complex matrix. Let $\rho \in \C^{m}, \gamma \in \C^n$ be two complex vectors
with all entries non-zero.
We denote by $$ \SC(A,\rho,\gamma)$$ the matrix obtained from $A$ by multiplying the $(i,j)$'th element of $A$ by $\rho_i \cdot \gamma_j$. We say that two matrices $A,B$ of the same dimensions are a scaling of each other if there exist non-zero vectors $\rho,\gamma$ such that $B = \SC(A,\rho,\gamma)$.
It is easy to check that this is an equivalence relation. We refer to the elements of the vector $\rho$ as the {\em row scaling coefficients} and to the elements of $\gamma$ as the {\em column scaling coefficients}. Notice that two matrices which are a scaling of each other have the same rank and the same pattern of zero and non-zero entries.
\end{define}

Matrix scaling originated in a paper of Sinkhorn \cite{Sinkhorn} and has been
widely studied since (see \cite{LSW98} for more background).
The following is a special case of a theorem from \cite{RS89} that gives sufficient conditions for finding a scaling of a matrix which has certain row and column sums.
\begin{define}[Property-$S$]
Let $A$ be an $m \times n$ matrix over some field. We say that $A$ satisfies {\em Property-$S$} if for every zero sub-matrix of $A$ of size $a \times b$ it holds that
\begin{equation}\label{eq-propS}
 \frac{a}{m} + \frac{b}{n} \leq 1.
\end{equation}
\end{define}

\begin{thm}[Matrix scaling theorem, Theorem 3 in \cite{RS89} ]\label{thm-scaling}
Let $A$ be an $m \times n$ real matrix with non-negative entries which satisfies Property-$S$. Then, for every $\eps >0$, there
exists a scaling $A'$ of $A$ such that the sum of each row of $A'$ is at most  $1+\eps$ and the sum of each column of $A'$ is at
least $m/n - \eps$. Moreover, the scaling coefficients used to obtain $A'$ are all positive real numbers.
\end{thm}

The proof of the theorem is algorithmic \cite{Sinkhorn}:
Start by normalizing $A$'s rows to have sum $1$,
then normalize $A$'s columns to have sum $m/n$,
then go back to normalizing the rows the have sum $1$, and so forth.
It can be shown (using a suitable potential function)
that this process eventually transforms $A$ to the claimed form
(since $A$ has Property-$S$).

We will use the following easy corollary of the above theorem.

\begin{cor}[$\ell_2^2$-scaling]\label{cor-l2scale}
Let $A = (a_{ij})$ be an $m \times n$ complex matrix which satisfies Property-$S$. Then, for every $\eps > 0$, there exists a scaling $A'$ of $A$ such that for every $i \in [m]$ $$ \sum_{j \in [n]} |a_{ij}|^2 \leq 1+\eps$$ and for every $j \in [n]$ $$ \sum_{i \in [m]} |a_{ij}|^2 \geq m/n - \eps.$$
\end{cor}
\begin{proof}
Let $B = (b_{ij}) = ( |a_{ij}|^2 )$. Then $B$ is a real non-negative matrix satisfying Property-$S$. Applying Theorem~\ref{thm-scaling} we get that for all $\eps>0$ there exists a scaling $B' = \SC(B, \rho, \gamma)$, with $\rho,\gamma$ positive real vectors,  which has row sums at most $1+\eps$ and column sums at least  $m/n - \eps$. Letting $\rho'_i = \sqrt{\rho_i}$ and $\gamma'_i = \sqrt{\gamma_i}$ we get a scaling $\SC(A,\rho',\gamma')$ of $A$ with the required properties.
\end{proof}

We will use a variant of a well known lemma (see for example \cite{Alo09}) which provides a bound on the rank of matrices whose diagonal entries are much larger than the off-diagonal ones.

\begin{lem}\label{lem-diagdom}
Let $A = (a_{ij})$ be an $n \times n$ complex hermitian matrix and let $0 < \ell < L$ be integers.
Suppose that $a_{ii} \geq L$ for all $i \in [n]$ and that $|a_{ij}| \leq \ell$ for all $i \neq j$.
Then
\[ \rank(A) \geq \frac{n}{1+ n\cdot (\ell/L)^2} \geq n - (n\ell/L)^2. \]
\end{lem}
\begin{proof}
We can assume w.l.o.g. that $a_{ii} = L$ for all $i$.
If not, then we can make the inequality into an equality by multiplying the $i$'th row and column by $(L/a_{ii})^{1/2} < 1$ without changing the rank or breaking the symmetry. Let $r = \rank(A)$ and let $\lambda_1,\ldots,\lambda_r$ denote the non-zero  eigenvalues of $A$ (counting multiplicities). Since $A$ is hermitian we have that the $\lambda_i$'s are real. We have
\begin{eqnarray*}
 n^2 \cdot L^2 &=& \tr(A)^2
 = \left( \sum_{i=1}^r \lambda_i \right)^2
 \leq  r \cdot \sum_{i=1}^r \lambda_i^2
 = r \cdot \sum_{i,j=1}^n |a_{ij}|^2 \\
 &\leq& r \cdot ( n\cdot L^2 + n^2 \cdot \ell^2).
\end{eqnarray*}
Rearranging we get the required bound. The second inequality in the statement of the lemma follows from the fact that $1/(1+x) \geq 1 - x$ for all $x$.
\end{proof}

\subsection{Proof of Theorem~\ref{thm-rankdesign}}
To prove the theorem we will first find a scaling of $A$ so that the norms (squared) of the columns are large and such that each entry is small.

Our first step is to find an $nk \times n$ matrix $B$ that will satisfy Property-$S$ and will be composed from rows of $A$ s.t. each row is repeated with multiplicity between $0$ and $q$. To achieve this  we will describe an algorithm that builds the matrix $B$ iteratively by concatenating to it rows from $A$. The algorithm will
\emph{mark} entries of $A$ as it continues to add rows. Keeping track of these marks will help us decide which rows to add next. Initially all the entries of $A$ are \emph{unmarked}. The algorithm proceeds in $k$ steps. At  step $i$ ($i$ goes from $1$ to $k$) the algorithm picks $n$ rows from $A$ and adds them to $B$. These $n$ rows are chosen as follows: For every $j \in \{1,\ldots,n\}$ pick a row that has an unmarked non-zero entry in the $j$'th column and mark this non-zero entry. The reason why such a row exists at all steps is that each column contains at least $k$ non-zero entries, and in each step we mark at most one non-zero entry in each column.
\begin{claim}\label{cla-matrixB}
The matrix $B$ obtained by the algorithm has Property-$S$ and each row of $A$ is added to $B$ at most $q$ times.
\end{claim}
\begin{proof}
The $n$ rows added at each of the $k$ steps form an $n\times n$ matrix with non-zero diagonal.
Thus they satisfy Property-$S$. It is an easy exercise to verify that a concatenation of matrices with Property-$S$ also has this property. The bound on the number of times each row is added to $B$ follows from the fact that each row has at most $q$ non-zero entries and each time we add a row to $B$ we mark one of its non-zero entries.
\end{proof}

Our next step is to obtain a scaling of $B$ and, from it, a scaling of $A$.
Fix some $\eps>0$ (which will later tend to zero).
Applying Corollary~\ref{cor-l2scale} we get a scaling $B'$ of $B$ such that the $\ell_2$-norm of each row is at most $\sqrt{1+\eps}$ and the $\ell_2$-norm of each column is at least $\sqrt{nk/n - \eps}= \sqrt{k -\eps}$.
We now obtain a scaling $A'$ of $A$ as follows:
The scaling of the columns are the same as for $B'$.
For the rows of $A$ appearing in $B$
we take the maximal scaling coefficient used for these rows in $B'$, that is,
if row $i$ in $A$ appears as rows $i_1,i_2,\ldots,i_{q'}$ in $B$,
then the scaling coefficient of row $i$ in $A'$ is the maximal scaling coefficient
of rows $i_1,i_2,\ldots,i_{q'}$ in $B'$.
For rows {\em not} in $B$,
we pick scaling coefficients so that their $\ell_2$ norm (in the final scaling) is equal to $1$.
\begin{claim}
The matrix $A'$ is a scaling of $A$ such that each row has $\ell_2$-norm at most
$\sqrt{1+\eps}$ and each column has $\ell_2$-norm at least $\sqrt{(k-\eps)/q}$.
\end{claim}
\begin{proof}
The fact that the row norms are at most $\sqrt{1+\eps}$ is trivial.
To argue about the column norms observe that a column of $B'$ is obtained from repeating each non-zero element in the corresponding column of $A'$ at most $q$ times (together with some zeros). Therefore, if we denote by $c_1,\ldots,c_s$ the non-zero entries in some column of $A'$, we have that
\[  \sum_{i=1}^s m_i \cdot |c_i|^2 \geq k - \eps, \]
where the $m_i$'s are integers between $0$ and $q$. In this last inequality we also relied on the fact that we chose the maximal row scaling coefficient among all those that correspond to the same row in $A$.  Therefore,
\[ \sum_{i=1}^s |c_i|^2 \geq (k-\eps)/q, \]
as required.
\end{proof}

Our final step is to argue about the rank of $A'$ (which is the same as the rank of $A$). To this end, consider the matrix $$ M = (A')^* \cdot A',$$
where $(A')^*$ is $A'$ transposed conjugate.
Then $M = (m_{ij})$ is an $n \times n$ hermitian matrix.
The diagonal entries of $M$ are exactly the squares of the $\ell_2$-norm of the columns of $A'$.
Therefore, $$ m_{ii} \geq  (k-\eps)/q$$ for all $i \in [n]$.

We now upper bound the off-diagonal entries.
The off-diagonal entries of $M$ are the inner products of different columns of $A'$.
The intersection of the support of each pair of different columns is at most $t$.
The norm of each row is at most $\sqrt{1+\eps}$.
For every two real numbers $\alpha,\beta$ so that $\alpha^2 + \beta^2 \leq 1+\eps$
we have $|\alpha \cdot \beta| \leq 1/2 + \eps'$, where $\eps'$ tends to zero as $\eps$ tends to zero.
Therefore
$$ |m_{ij}| \leq t\cdot (1/2+\eps')$$ for all $i \neq j \in [n]$. Applying Lemma~\ref{lem-diagdom} we get that
$$\rank(A) = \rank(A') \geq n - \left(\frac{q\cdot t(1/2+\eps') \cdot n}{k-\eps}\right)^2.$$
Since this holds for all $\eps > 0$
it holds also for $\eps = 0$, which gives the required bound on the rank of $A$. \qed

\section{Sylvester-Gallai Configurations}\label{sec-SG}

In this section we prove the quantitative Sylvester-Gallai (SG) Theorem. We will be interested with point configurations in real
and complex space. These are finite sets of distinct points $v_1,\ldots,v_n$ in $\R^d$ or $\C^d$. The dimension of a
configuration is defined to be the dimension of the smallest affine subspace containing all points.

\begin{define}[Special and ordinary lines]
Let $v_1,\ldots,v_n \in \C^d$ be a set of $n$ distinct points in $d$-dimensional complex space. A line $\ell$ passing through at
least three of these points is called a {\em special} line. A line passing through exactly two points is called an {\em ordinary}
line.
\end{define}

\begin{define}[$\delta$-SG configuration]
\label{def: d SG} Let $\delta \in [0,1]$. A set of $n$ distinct points $v_1,\ldots,v_n \in \C^d$ is called a {\em $\delta$-SG
configuration} if for every $i \in [n]$, there exists a family of special lines $L_i$ all passing through $v_i$ and at least
$\delta n$ of the points $v_1,\ldots,v_n$ are on the lines in $L_i$. (Note that each collection $L_i$ may cover a different
subset of the $n$ points.)
\end{define}

The main result of this section bounds the dimension of $\delta$-SG configurations for all $\delta > 0$. Since we can always
satisfy the definition by spreading the points evenly over  $1/\delta$ lines we know that the dimension can be at least
$2/\delta$ (and in fact in complex space at least $3/\delta$). We prove an upper bound of $O(1/\delta^2)$.

\begin{thm}[Restatement of Theorem~\ref{ithm:SG} -- quantitative SG theorem]\label{thm-deltaSG}
Let $\delta \in (0,1]$. Let $v_1,\ldots,v_n \in \C^d$ be a $\delta$-SG configuration. Then
$$\dimension\{v_1,\ldots,v_n\} < 13 / \delta^2.$$
Moreover, the dimension of a $1$-SG configuration is at most $10$.
\end{thm}

The constants in the proof have been optimized to the best of our abilities. Notice that in the above theorem $\delta$ can be
dependant on $n$. For example,  a $ (1/\log(n))$-SG configuration of $n$ points can have rank at most $O( \log(n)^2 )$.

\subsection{Preliminaries to the Proof of Theorem~\ref{thm-deltaSG}}

The notion of a latin square will turn out useful in the proof:

\begin{define}[Latin squares]
An $r \times r$ {\em latin square} is an $r\times r$ matrix $D$ such that $D_{i,j}\in [r]$ for all $i,j$ and every number in
$[r]$ appears exactly once in each row and in each column. A latin square $D$ is called {\em diagonal} if $D_{i,i}=i$ for all
$i\in[r]$.
\end{define}

\begin{thm}[\cite{Hil73}]\label{thm-diagonalsquare}
For every $r \geq 3$ there exists a diagonal $r \times r$ latin square.
\end{thm}

We note that we use diagonal latin squares only to optimize constant factors.
If one does not care about such factors
then there is a simple construction that serves the same goal.


The following lemma is an easy consequence of the above theorem.

\begin{lem}\label{lem-triples}
Let $r \geq 3$. 
Then there exists a set  $T \subset [r]^3$ of $r^2 - r$ triples that satisfies the following properties:
\begin{enumerate}
\item Each triple $(t_1,t_2,t_3) \in T$ is of three distinct elements.
\item For each $i \in [r]$ there are exactly $3(r-1)$ triples in $T$ containing $i$ as an element.
\item For every pair $i,j \in [r]$ of distinct elements there are at most $6$ triples in $T$ which contain both $i$ and $j$ as elements.
\end{enumerate}
\end{lem}
\begin{proof}
Let $D$ be an $r \times r$ diagonal latin square which we know exists from Theorem~\ref{thm-diagonalsquare}. Define $T \subset
[r]^3$ to be the set of all triples $(i,j,k) \in [r]^3$ with $i \neq j$ such that $D_{i,j}=k$. The number of such triples is $r^2
- r$. Property $1$ holds by the definition of diagonal latin square--- we cannot have $D_{i,j}=i$ for $j\neq i$ since $D_{i,i}=i$
and every row in $D$ has distinct as the $(i,i)$ entry in $D$ is labeled $i$ for all $i \in [r]$, and similarly we cannot have
$D_{i,j}=j$ for $i\neq j$.

Let $i\in [r]$. By construction, there are $r-1$ triples in $T$ which have $i$ as their first entry, and $r-1$ triples that have
$i$ as their second entry. There are also $r-1$ triples in $T$ which have $i$ as their last entry, since for every one of the
$r-1$ rows $i'\neq i$ there is exactly one location $j'\neq i'$ in which the label $i$ appears, and that contributes the triple
$(i',j',i)$ to $T$. This proves Property~$2$.

To prove Property~$3$ observe that two triples in $T$ can agree in at most one place. For example, knowing the row and column
determines the label, knowing the row and label determines the column, and so forth. Therefore, a pair $(i,j)$ cannot appear in
more than $6$ triples since otherwise there would have been at least two triples with $i,j$ at the same places, and these triples
would violate the above rule.
\end{proof}

\subsection{Proof of Theorem~\ref{thm-deltaSG}}

Let $V$ be the $n \times d$ matrix whose $i$'th row is the vector $v_i$. Assume w.l.o.g. that $v_1 = 0$. Thus
$$\dimension\{v_1,\ldots,v_n\} = \rank(V) .$$

The overview of the proof is as follows. We will first build an $m \times n$ matrix $A$ that will satisfy $A \cdot V = 0$. Then,
we will argue that the rank of $A$ is large because it is a design matrix. This will show that the rank of $V$ is small.

Consider a special line $\ell$ which passes through three points $v_i,v_j,v_k$. This gives a linear dependency among the three
vectors $v_i,v_j,v_k$ (we identify a point with its vector of coordinates in the standard basis). In other words, this gives a
vector $a = (a_1,\ldots,a_n)$ which is non-zero only in the three coordinates $i,j,k$ and such that $a \cdot V = 0$. If $a$ is
not unique, choose an arbitrary vector $a$ with these properties.

Our strategy is to pick a family of collinear triples among the points in our configuration and to build the matrix $A$ from rows
corresponding to these triples in the above manner.

Let $\L$ denote the set of all special lines in the configuration (i.e. all lines containing at least three points). Then each
$L_i$ is a subset of $\L$ containing lines passing through $v_i$.  For each $\ell \in \L$ let $V_\ell$ denote the set of points
in the configuration which lie on the line $\ell$. Then $|V_\ell| \geq 3$ and we can assign to it a family of triples $T_\ell
\subset V_\ell^3$, given by Lemma~\ref{lem-triples} (we identify $V_\ell$ with $[r]$, where $r = |V_\ell|$ in some arbitrary
way).

We now construct the matrix $A$ by going over all lines $\ell \in \L$ and for each triple in $T_\ell$ adding as a row of $A$ the
vector with three non-zero coefficients $a = (a_1,\ldots,a_n)$ described above (so that $a$ is the linear dependency between the
three points in the triple).

Since the matrix $A$ satisfies $A \cdot V = 0$ by construction, we only have to argue that $A$ is a design matrix and bound its
rank.
\begin{claim}
The matrix $A$ is a $(3,3k,6)$-design matrix,  where $k \triangleq \lfloor \delta n \rfloor - 1$.
\end{claim}
\begin{proof}

By construction, each row of $A$ has exactly $3$ non-zero entries. The number of non-zero entries in column $i$ of $A$
corresponds to the number of triples we used that contain the point $v_i$. These can come from all special lines containing
$v_i$. Suppose there are $s$ special lines containing $v_i$ and let $r_1,\ldots,r_s$ denote the number of points on each of those
lines. Then, since the lines through $v_i$ have only the point $v_i$ in common, we have that
\[ \sum_{j=1}^s (r_j - 1) \geq k. \]
The properties of the families of triples $T_\ell$ guarantee that there are $3(r_j-1)$ triples containing $v_i$ coming from the $j$'th line. Therefore there are at least $3k$ triples in total containing $v_i$.

The size of the  intersection of columns $i_1$ and $i_2$ is equal to the number of triples containing the points $v_{i_1},v_{i_2}$ that were used in the construction of $A$. These triples can only come from one special line (the line containing these two points) and so, by Lemma~\ref{lem-triples}, there can be at most $6$ of those.
\end{proof}

Applying Theorem~\ref{thm-rankdesign} we get that
\begin{eqnarray*}
 \rank(A) & \geq & n - \left(\frac{3 \cdot 6 \cdot n}{2 \cdot 3k}\right)^2
\geq n - \left(\frac{3 \cdot n}{\delta n  - 2}\right)^2 \\
 & \geq & n - \left(\frac{3 \cdot n \cdot 13 }{ 11 \cdot \delta n} \right)^2 >   n - 13/\delta^2,
\end{eqnarray*}
where the third inequality holds as $\delta n \geq 13$
since otherwise the theorem trivially holds.
Since $A \cdot V= 0$ we have that
\[ \rank(A) + \rank(V) \leq n .\]
This implies that
\[ \rank(V) < 13/\delta^2 , \]
which completes the proof.
For $\delta = 1$, the calculation above yields $\rank(V) < 11$.
\qed

\subsection{Average-Case Version}

In this section we use Theorem~\ref{thm-deltaSG} to argue about the case where we only know that there are many collinear triples in a configuration.

\begin{thm}[Average-case SG theorem]\label{thm-SGav}
Let $V = \{v_1,\ldots,v_m\} \subset \C^d$ be a set of $m$ distinct points. Let $T$ be the set of (unordered)
collinear triples in $V$. Suppose $|T| \geq \alpha m^2$
and that every two points $v,v'$ in $V$ appear in at most $c$ triples in $T$,
then there exists a subset $V' \subset V$ such that
$|V'| \geq \alpha m/(2c)$ and $\dimension(V') \leq O(1/\alpha^2)$.
\end{thm}
Notice that the bound on the number of triples containing a fixed pair of points is necessary for the theorem to hold. If we remove this assumption than we could create a counter-example by arranging the points so that $m^{2/3}$ of them are on a line and the rest span the entire space.

\begin{lem}
\label{lem: comb lem}
Let $H$ be a $3$-regular hypergraph with vertex set $[m]$ and $\alpha m^2$ edges
of co-degree at most $c$ (i.e. for every $i \neq j$ in $[m]$,
the set $\{i,j\}$ is contained in at most $c$ edges).
Then there is a subset $M \subseteq [m]$ of size $|M| \geq \alpha m/(2c)$
so that the minimal degree of the sub-graph of $H$ induced by $M$
is at least $\alpha m / 2$.
\end{lem}

\begin{proof}
We describe an iterative process to find $M$.
We start with $M = [m]$.
While there exists a vertex of degree less than $\alpha m / 2$,
remove this vertex from $M$ and remove all edges containing this vertex from $H$.
Continuing in this fashion we conclude with a set $M$ such that every point in
$M$ has degree at least $\alpha m / 2$.
This process removed in total at most $m \cdot \alpha m / 2$ edges and thus the new
$H$ still contains at least $\alpha m^2 / 2$ edges.
As the co-degree is at most $c$,
every vertex appears in at most $cm$ edges.
Thus, the size of $M$ is of size at least $\alpha m / (2c)$.
\end{proof}

\begin{proof}[Proof of Theorem~\ref{thm-SGav}]
The family of triples $T$ defines a $3$-regular hypergraph on $V$
of co-degree at most $c$.
Lemma~\ref{lem: comb lem} thus implies that
there is a subset $V' \subseteq V$ of size $|V'| \geq \alpha m/(2c)$
that is an $(\alpha/2)$-SG configuration. By Theorem~\ref{thm-deltaSG},
$V'$ has dimension at most $O(1/\alpha^2)$.
\end{proof}

%
%
%
%
%
%

\section{Robust SG Theorem for $k$-Flats}\label{sec-highdim}

In this section we prove two high-dimensional analogs of the SG theorem.
Let $\fl(v_1,\ldots,v_k)$ (fl for `flat') denote the affine span of $k$ points (i.e. the points that can be written as linear combinations with coefficients that sum to one). We call $v_1,\ldots,v_k$ \emph{independent} if their flat is of dimension $k-1$ (dimension means affine dimension),
and say that $v_1,\ldots,v_k$ are \emph{dependent} otherwise.
A \emph{$k$-flat} is an affine subspace of dimension $k$.

In the following $V$ is a set of $n$ distinct points in complex space $\C^d$.
A $k$-flat is called \emph{ordinary} if its intersection with $V$ is contained in the union of a $(k-1)$-flat and a single point. A $k$-flat is
\emph{elementary} if its intersection with $V$ has exactly $k+1$ points. Notice that for $k=1$ (lines) the two notions of ordinary and elementary coincide.

For dimensions higher than one,
there are two different definitions that generalize that of $\SG$ configuration.
The first definition is based on ordinary $k$-flats
(though in a slightly stronger way which will be more useful in the proofs to come).
The second definition (which is less restricted than the first one) 
uses elementary $k$-flats.

\begin{define}
The set $V$ is a $\delta$-$\SG_k^*$ configuration if for every independent $v_1,\ldots,v_k \in V$ there are at least $\delta n$ points $u \in V$ s.t. either $u \in \fl(v_1,\ldots,v_k)$ or the $k$-flat $\fl(v_1,\ldots,v_k,u)$ contains a point $w$ outside $\fl(v_1,\ldots,v_k) \cup \{ u\}$.
\end{define}

\begin{define}
The set $V$ is a $\delta$-$\SG_k$ configuration if for every independent $v_1,\ldots,v_k \in V$ there are at least $\delta n$ points $u \in V$ s.t. either $u \in \fl(v_1,\ldots,v_k)$ or the $k$-flat $\fl(v_1,\ldots,v_k,u)$ is not elementary.
\end{define}

Both definitions coincide with that of $\SG$ configuration when $k = 1$:
Indeed, $\fl(v_1) = v_1$ and $\fl(v_1,u)$ is the line through $v_1,u$.
Therefore, $u$ is never in $\fl(v_1)$ and the line $\fl(v_1,u)$
is not elementary iff it contains at least one point $w \not \in \{v_1,u\}$.

We prove two high-dimensional versions of the SG theorem,
each corresponding to one of the definitions above.
The first uses the more restricted `star' definition and
gives a strong upper bound on dimension.
The second uses the less restricted definition
and gives a weaker bound on dimension.

\begin{thm}\label{thm-weak}
Let $V$ be a $\delta$-$\SG_k^*$ configuration. Then $dim(V) \leq f(\delta,k)$ with
$$ f(\delta,k) = O\left( (k/\delta)^{2} \right).$$
\end{thm}

\begin{thm}
\label{thm-strong}
Let $V$ be a $\delta$-$\SG_k$ configuration. Then $dim(V) \leq g(\delta,k)$ with
$$g(\delta,k) = 2^{C^k} / \delta^2$$
with $C > 1$ a universal constant.
\end{thm}

The proofs of the two theorems are below.
Theorem~\ref{thm-weak} follows by an appropriate induction
on the dimension, using the (one-dimensional) robust SG theorem.
Theorem~\ref{thm-strong} follows by reduction
to Theorem~\ref{thm-weak}.

Before proving the theorems we set some notations.
Fix some point $v_0 \in V$. By a \emph{normalization w.r.t. $v_0$} we mean an affine transformation $N : \C^d \mapsto \C^d$ which first moves $v_0$ to zero, then picks a hyperplane $H$ s.t. no point in $V$ (after the shift) is parallel to $H$ (i.e has inner product zero with the orthogonal vector to $H$) and finally multiplies each point (other than zero) by a constant s.t. it is in $H$.

\begin{claim}
\label{clm: prop of N}
For such a mapping $N$ we have that  $v_0,v_1,\ldots,v_k$ are dependent
iff $N(v_1),\ldots,N(v_k)$ are dependent.
\end{claim}

\begin{proof}
Since translation and scaling does not affect dependence,
w.l.o.g. we assume that $v_0 = 0$ and that
the distance of the hyperplane $H$ from zero is one.
Let $h$ be the unit vector orthogonal to $H$.
For all $i \in [k]$ we have $N(v_i) = v_i / \ip{v_i}{h}$.
Assume that $v_0,v_1,\ldots,v_k$ are dependent, that is,
w.l.o.g. $v_k = \sum_{i \in [k-1]} a_i v_i$ for some $a_1,\ldots,a_{k-1}$.
For all $i \in [k-1]$ define $b_i = a_i \ip{v_i}{h} / \ip{v_k}{h}$.
Thus $N(v_k) = \sum_{i \in [k-1]} a_i v_i / \ip{v_k}{h} =
\sum_{i \in [k-1]} b_i N(v_i)$
where $\sum_{i \in [k-1]} b_i = 1$,
which means that $N(v_1),\ldots,N(v_k)$ are dependent.
Since the map $a_i \mapsto b_i$ is invertible,
the other direction of the claim holds as well.
\end{proof}

We first prove the theorem for $\delta$-$\SG_k^*$ configurations.

\begin{proof}[Proof of Theorem~\ref{thm-weak}]
The proof is by induction on $k$. For $k=1$ we know $f(\delta,1)
\leq c \delta^{-2}$ with $c > 1$ a universal constant.
Suppose $k > 1$. We separate into two cases. The first case is when $V$ is an $(\delta/(2k))$-$\SG_1$ configuration and we are done using the bound on $k=1$. In the other case there is some point $v_0 \in V$ s.t. the size of the set of points on special lines through $v_0$ is at most $\delta/(2k)$ (a line is special if it
contains at least three points). Let $S$ denote the set of points on special lines through $v_0$. Thus $|S| < \delta n / (2k)$. Let $N : \C^d \mapsto \C^d$ be a normalization w.r.t. $v_0$. Notice that for points $v \not\in S$ the image $N(v)$ determines $v$. Similarly, all points on some special line map to the same point via $N$.

Our goal is to show that $V' = N(V \setminus \{v_0\})$ is a
$((1-1/(2k))\delta)$-$SG_{k-1}^*$ configuration
(after eliminating multiplicities from $V'$).
This will complete the proof since $\dim(V) \leq \dim(V') + 1$.
Indeed, if this is the case we have
$$f(\delta,k) \leq \max \{ 4c (k/\delta)^2 , f((1-1/(2k))\delta,k-1) + 1\}.$$
and by induction we have $f(\delta,k) \leq 4 c (k/\delta)^2$.

Fix $v'_1,\ldots,v'_{k-1} \in V'$ to be $k-1$ independent points
(if no such tuple exists then $V'$ is trivially a configuration).
Let $v_1,\ldots,v_{k-1} \in V$ be points s.t. $N(v_i) = v'_i$ for $i \in [k-1]$.
Claim~\ref{clm: prop of N} implies that $v_0,v_1,\ldots,v_{k-1}$ are independent. Thus, there is a set $U \subset V$ of size at least $\delta n$ s.t.
for every $u \in U$ either $u \in \fl(v_0,v_1,\ldots,v_{k-1})$ or the $k$-flat $\fl(v_0,v_1,\ldots,v_{k-1},u)$ contains a point $w$ outside $\fl(v_0,v_1,\ldots,v_{k-1}) \cup \{ u \}$.

Let $\tilde U = U \setminus S$ so that $N$ is invertible on $\tilde U$
and $$|\tilde U| \geq |U| - |S| \geq (1-1/(2k)) \delta n.$$
Suppose $u \in \tilde U$ and let $u' = N(u)$.
By Claim~\ref{clm: prop of N} if $u \in \fl(v_0,v_1,\ldots,v_{k-1})$ then $u'$ is in
$\fl(v'_1,\ldots,v'_{k-1})$.
Otherwise, $\fl(v_0,v_1,\ldots,v_{k-1},u)$ contains a point $w$ outside $\fl(v_0,v_1,\ldots,v_{k-1}) \cup \{ u \}$. Let $w' = N(w)$.
We will show that $w'$ is (a) contained in the $(k-1)$-flat $\fl(v'_1,\ldots,v'_{k-1},u')$ and (b) is outside $\fl(v'_1,\ldots,v'_{k-1}) \cup \{ u' \}$.
Property (a) follows from Claim~\ref{clm: prop of N} since $v_0$, $v_1,\ldots,v_{k-1}$, $u,w$ are dependent and so $v'_1,\ldots,v'_{k-1},u',w'$ are also dependent. To show (b) observe first that by Claim~\ref{clm: prop of N}
the points $v'_1,\ldots,v'_{k-1},u'$ are independent (since $v_0,v_1,\ldots,v_{k-1},u$ are independent) and so $u'$ is not in $\fl(v'_1,\ldots,v'_{k-1})$. We also need to show that $w' \neq u'$ but this follows from the fact that $u \neq w$
and so $w' = N(w) \neq N(u) = u'$ since $N$ is invertible on $\tilde U$
and $u \in \tilde U$.
Since $$|N(\tilde U)| = |\tilde U| \geq (1-1/(2k)) \delta n \geq (1-1/(2k)) \delta |V'|$$ the proof is complete.
\end{proof}

We can now prove the theorem for $\delta$-$\SG_k$ configurations.

\begin{proof}[Proof of Theorem~\ref{thm-strong}]
The proof follows by induction on $k$ (the case $k=1$ is given by Theorem~\ref{thm-deltaSG}).
Suppose $k>1$. Suppose that $\dim(V) > g(\delta,k)$. We want to show that there exist $k$ independent points $v_1,\ldots,v_k$ s.t. for at least $1-\delta$ fraction of the points $w \in V$ we have
that $w$ is not in $\fl(v_1,\ldots,v_k)$ {\bf and} the flat $\fl(v_1,\ldots,v_k,w)$ is elementary (i.e. does not contain any other point).

Let $k' = g(1,k-1)$.
By choice of $g$ we have
$g(\delta,k) > f(\delta,k'+1)$
with $f$ from Theorem~\ref{thm-weak}.
Thus, by Theorem~\ref{thm-weak}, we can find $k'+1$ independent points $v_1,\ldots,v_{k'+1}$ s.t. there is a set $U \subset V$ of size at least $(1-\delta) n$ s.t. for every $u \in U$ we have that $u$ is not in $\fl(v_1,\ldots,v_{k'+1})$ {\bf and} the $(k'+1)$-flat $\fl(v_1,\ldots,v_{k'+1},u)$ contains only one point,
namely $u$, outside $\fl(v_1,\ldots,v_{k'+1})$.

We now apply the inductive hypothesis on the set $V \cap \fl(v_1,\ldots,v_{k'+1})$ which has dimension at least $k' = g(1,k-1)$.
This gives us $k$ independent points $v'_1,\ldots,v'_{k}$ that define an elementary $(k-1)$-flat $\fl(v'_1,\ldots,v'_k)$. (Saying that $V$ is not $1$-$\SG_{k-1}$ is the same as saying that it contains an elementary $(k-1)$-flat). Joining any of the points $u \in U$ to $v'_1,\ldots,v'_k$ gives us an elementary $k$-flat and so the theorem is proved.
\end{proof}

\section{Generalizations of the Motzkin-Rabin Theorem} 
\label{sec-MR}
In this section we prove two variants of the Motzkin-Rabin Theorem. The first is a quantitative analog in the spirit of Theorem~\ref{thm-deltaSG}. The second is a variant in which the number of colors is three (instead of two).

\subsection{A Quantitative Variant}

\begin{define}[$\delta$-MR configuration]
Let $V_1,V_2$ be two disjoint finite subsets of $\C^d$.
Points in $V_1$ are of \emph{color} $1$
and points in $V_2$ are of \emph{color} $2$.
A line is called \emph{bi-chromatic} if it contains
at least one point from each of the two colors.
We say that $V_1,V_2$ are a {\em $\delta$-MR} configuration if
for every $i \in [2]$ and for every point $p \in V_i$,
the bi-chromatic lines through $p$ contain at least $\delta |V_i|$ points.
\end{define}

\begin{thm}
Let $V_1,V_2 \subset \C^d$ be a $\delta$-MR configuration. Then $$\dimension(V_1,V_2) \leq O(1/\delta^4).$$
\end{thm}
\begin{proof}
We will call a line passing through exactly two points in $V_1$ (resp. $V_2$) a $V_1$-ordinary (resp. $V_2$-ordinary) line.
W.l.o.g. assume $$|V_1| \leq |V_2|.$$ We seperate the proof into two cases:

Case I is when $V_2$ is a $(\delta/2)$-SG configuration. Then, by Theorem~\ref{thm-deltaSG}, $\dimension(V_2) \leq
O(1/\delta^2)$. If in addition  $$\dimension(V_1) \leq 13/(\delta/2)^2$$ then we are done. Otherwise, by
Theorem~\ref{thm-deltaSG}, there exists a point $a_0 \in V_1$ such that there are at least $(1- \delta/2)|V_1|\,\,\,$
$V_1$-ordinary lines through $a_0$. Let $a_1,\ldots,a_k$ denote the points in $V_1$ that belong to these lines with $k \geq
(1-\delta/2)|V_1|$. We now claim that $V_2 \cup \{a_0\}$ spans all the points in $V_1$. This will suffice since, in this case,
$\dimension(V_2) \leq O(1/\delta^2)$. Let $a \in V_1$. Then, since $V_1,V_2$ is a $\delta$-MR configuration, there are at least
$\delta |V_1|$ points in $V_1$ such that the line through them and $a$ contains a point in $V_2$. One of these points must be
among $a_1,\ldots,a_k$, say it is $a_1$. Since $a$ is in the span of $V_2$ and $a_1$ and since $a_1$ is in the span of $V_2$ and
$a_0$ we are done.

Case II is when $V_2$ is not a $(\delta/2)$-SG configuration. In this case, there is a point $b \in V_2$ such that there are at
least $(1 - \delta/2)|V_2|\,\,$ $V_2$-ordinary lines through $b$. From this fact and from the $\delta$-MR property, we get that
$|V_1| \geq (\delta/2)|V_2|$ (there are at least $(\delta/2)|V_2|\,\,$ $V_2$-ordinary lines through $b$ that have an additional
point from $V_1$ on them). This implies that the union $V_1 \cup V_2$ is a $(\delta^2/4)$-SG configuration and the result follows
by applying Theorem~\ref{thm-deltaSG}.
\end{proof}

\subsection{A Three Colors Variant}

\begin{define}[3MR configuration]
Let $V_1,V_2,V_3$ be three pairwise disjoint finite subsets of $\C^d$,
each of distinct points.
We say that $V_1,V_2,V_3$ is a {\em 3MR}-configuration if
every line $\ell$ so that $\ell \cap (V_1\cup V_2 \cup V_3)$ has more than one point
intersects at least two of the sets $V_1,V_2,V_3$.
\end{define}

\begin{thm}
Let $V_1,V_2,V_3$ be a 3MR configuration
and denote $V = V_1 \cup V_2 \cup V_3$.
Then $$\dimension(V) \leq O(1).$$
\end{thm}
\begin{proof}
Assume w.l.o.g. that $V_1$ is not smaller than $V_2,V_3$.
Let $\alpha = 1/16$.
There are several cases to consider:

\begin{description}
\item[1. $V_1$ is an $\alpha$-SG configuration.]
By Theorem~\ref{thm-deltaSG},
the dimension of $V_1$ is at most
$$d_1 = O(1/\alpha^2).$$
Consider the two sets
$$V'_2 = V_2 \setminus \span(V_1) \ \ \text{and} \ \ V'_3 = V_3 \setminus \span(V_1),$$
each is a set of distinct points in $\C^{d}$.
Assume w.l.o.g. that $|V'_2| \geq |V'_3|$.

\begin{description}
\item[1.1. $V'_2$ is an $\alpha$-SG configuration.]
By Theorem~\ref{thm-deltaSG},
the dimension of $V'_2$ is at most
$$d_2 = O(1/\alpha^2).$$
Fix a point $v_3$ in $V'_3$.
For every point $v \neq v_3$ in $V'_3$
the line through $v_3,v$ contains a point from $\span(V_1) \cup V'_2$.
Therefore,
$$\dimension(V) \leq d_1 + d_2 + 1 \leq O(1).$$

\item[1.2. $V'_2$ is not an $\alpha$-SG configuration.]
There is a point $v_2$ in $V'_2$ so that
for $k \geq |V'_2| / 2$ of the points $v \neq v_2$ in $V'_2$
the line through $v_2,v$ does not contain any other point from $V'_2$.
If $V'_2 = \span(V_1,v_2)$ then the dimension of $V_1 \cup V_2$ is at most
$d_1 + 1$ and we are done as in the previous case.
Otherwise, there is a point $v'_2$ in $V'_2 \setminus \span(V_1,v_2)$.

We claim that in this case $|V'_3| \geq k /2$.
Denote by $P_2$ the $k$ points $v \neq v_2$ in $V'_2$
so that the line through $v_2,v$ does not contain any other point from $V'_2$.
For every $v \in P_2$ there is a point $V_{1,3}(v)$ in $V_1 \cup V_3$
that is on the line through $v,v_2$ (the point $v_2$ is fixed).
There are two cases to consider.

The first case is that for at least $k/2$ of the points $v$ in $P_2$
we have $V_{1,3}(v) \in V_3$.
In this case clearly $|V_3| \geq k/2$.

The second case is that for at least $k/2$ of the points $v$ in $P_2$
we have $V_{1,3}(v) \in V_1$.
Fix such a point $v \in P_2$ (which is in $\span(V_1,v_2)$).
The line through $v'_2,v$ contains a point $v'$ from $V_1 \cup V_3$.
The point $v'$ is not in $\span(V_1)$,
as if it was then $v'_2$ would be in $\span(v,v') \subseteq \span(V_1,v)$.
Therefore $v'$ is in $V_3$.
This also implies that $|V'_3| \geq k/2$.

Denote $V' = V_2 \cup V_3'$.
So we can conclude that
for every $v'$ in $V'$ the special lines through $v'$
contain at least $|V'|/8$ of the points in $V_1 \cup V_2 \cup V_3$.
As in the proof of Theorem~\ref{thm-deltaSG},
we can thus define a family of triples $T$, each triple of three distinct
collinear points in $V$, so that each $v'$ in $V'$ belongs to at least $|V'|/8$
triples in $T$ and each two distinct $v',v''$ in $V'$ belong to at most $6$ triples.

By a slight abuse of notation,
we also denote by $V$ the matrix with rows defined by the points in $V$.
Let $V_1$ be the submatrix of $V$ with row defined by points in $\span(V_1) \cap V$
and $V'$ be the submatrix of $V$ with row defined by points in $V'$.
Use the triples in $T$ to construct a matrix $A$ so that $A \cdot V = 0$.
Let $A_1$ be the submatrix of $A$ consisting of the columns that correspond to $\span(V_1) \cap V$ and
$A'$ be the submatrix of $A$ consisting of the columns that correspond to $V'$.
Therefore, $A' \cdot V' = - A_1 \cdot V_1$ which implies
$$\rank(A' \cdot V') \leq \rank(A_1 \cdot V_1) \leq d_1 .$$
By the above discussion $A'$ is a $(3, |V'|/8, 6)$-design matrix and thus, by Theorem~\ref{thm-rankdesign}, has rank at least $$ |V'| - O(1)$$ and so
$$\dimension(V') \leq O(1) + d_1 \leq O(1).$$
We can finally conclude that
$$\dimension(V)
\leq d_1 + \dimension(V') \leq O(1).$$

\end{description}

\item[2. $V_1$ is not an $\alpha$-SG configuration.]
There is a point $v_1$ in $V_1$ so that
for at least $|V_1| / 2$ of the points $v \neq v_1$ in $V_1$
the line through $v_1,v$ does not contain any other point from $V_1$.
Assume w.l.o.g. that $|V_2| \geq |V_3|$.
This implies that
$$|V_2| \geq |V_1| / 4 .$$

\begin{description}
\item[2.1. $|V_3| < |V_2| / 16$.]
In this case
the configuration defined by $V_1 \cup V_2$ is an $\alpha$-SG configuration.
By Theorem~\ref{thm-deltaSG},
the dimension of $V_1 \cup V_2$ is at most
$$d_{1,2} = O(1/\alpha^2).$$
Fix a point $v_3$ in $V_3$.
For every point $v \neq v_3$ in $V_3$
the line through $v_3,v$ contains a point from $V_1 \cup V_2$.
Therefore,
$$\dimension(V)
\leq d_{1,2} + 1 \leq O(1).$$

\item[2.1. $|V_3| \geq |V_2| / 16$.]
In this case $V$
is an $\alpha$-SG configuration.
By Theorem~\ref{thm-deltaSG},
the dimension of $V$ is thus at most
$O(1/\alpha^2)$.

\end{description}

\end{description}
\end{proof}

\section{Two-Query Locally Correctable Codes} 
\label{sec-LCC}

We now prove the non-existence of $2$-query (linear) locally correctable codes (LCC) over $\C$. We start by formally defining
locally correctable codes:

\begin{define}[Linear locally correctable code (LCC)]\label{def-lcc}
Let $\F$ be some field. A $(q,\delta)$-LCC over $\F$ is a linear subspace $C \subset \F^m$ such that there exists a randomized
decoding procedure $D : \F^m \times [m] \mapsto \F$ with the following properties:
\begin{enumerate}
\item For all $x \in C$, for all $i \in [m]$ and for all $v \in \F^m$ with $w(v) \leq \delta m$ we have that $D\left( x + v, i\right) = x_i$ with probability at least $3/4$ (the probability is taken only over the internal randomness of $D$).
\item For every $y \in \F^m$ and $i \in [m]$, the decoder $D(y,i)$ reads at most $q$ positions in $y$.
\end{enumerate}
The {\em dimension} of an LCC is simply its dimension as a subspace of $\F^m$.
\end{define}

In the above definition we allow the algorithm $D$ to perform operations over the  field $\F$. Since we do not care about the
running time of $D$ we do not discuss issues of representation of field elements and efficiency of handling them. (In any case,
it turns out that for linear codes in the small number of queries and low error case, one can assume w.l.o.g.
that the decoder is also linear, see Lemma~\ref{lem-LCC} below.)

Our result on locally decodable codes is the following:

\begin{thm}[Restatement of Theorem~\ref{ithm:lcc}--- non-existence of $2$ query LCCs over $\C$]\label{thm-lcc}
Let $C \subset \C^m$ be a $(2,\delta)$-LCC over $\C$. Then
\[ \dim(C) \leq O(1/\delta^9). \]
\end{thm}

As in Theorem~\ref{thm-deltaSG}, also in this theorem, $\delta$ can be an arbitrary function of $m$. To make the connection
between LCCs and $SG$-configurations explicit, we define the notion of  a $\delta$-LCC configuration.

\begin{define}[$\delta$-LCC Configuration]
A list of non-zero points $(v_1,\ldots, v_m)$ in $\C^d$ (not necessarily distinct) is called a $\delta$-LCC configuration if for
every subset $\Delta \subset [m]$ of size at most $\delta m$ and for every $i \in [m]$, there exist $j,k \in [m] \setminus
\Delta$ such that either $v_i \in \{v_j,v_k\}$ (in which case $v_i$ can be recovered by its own copies), or $v_i,v_j,v_k$ are
three distinct collinear points (in which case $v_i$ is recovered by two other coordinates).
\end{define}

The following lemma shows the connection between these two notions.

\begin{lem}\label{lem-LCC}
If there exists a $(2,\delta)$-LCC  
of dimension $n$ over $\C$ then there exists a $\delta$-LCC configuration
of dimension at least $n-1$ over $\C$.
\end{lem}

To prove the lemma we will use the following definition.

\begin{define}[Generating set]
Let $C \subset \F^m$ be a subspace. We say that a list of vectors $V = (v_1,\ldots,v_m)$ in $\F^n$ is a {\em generating set} for $C$ if
\[ C = \left\{ \left(\ip{y}{v_1}, \ip{y}{v_2}, \ldots, \ip{y}{v_m} \right) \,\,|\,\, y \in \F^n \right\} , \]
where $\ip{y}{v}$ is the standard inner product over $\F$.
\end{define}

\begin{proof}[Proof of Lemma~\ref{lem-LCC}]
Let $V = (v_1,\ldots,v_m)$ be a generating set for $C$ with $\dimension(V) \geq n-1$. We might lose $1$ since we defined
$\dimension(V)$ as the dimension of the smallest {\em affine} subspace containing $V$.  When the local decoder for $C$ reads two
positions in a codeword, it is actually reading $\ip{y}{v_j}, \ip{y}{v_k}$ for some vector $y \in \C^n$ (or noisy versions of
them). In order to be able to recover $\ip{y}{v_i}$ from $\ip{y}{v_j}, \ip{y}{v_k}$ with positive probability it must be that
$v_i \in \span\{v_j,v_k\}$.
(If we choose $y$ as Gaussian and $v_i$ is not in the span of $v_j,v_k$
then even conditioned on the values of $\ip{y}{v_j},\ip{y}{v_k}$
the r.v. $\ip{y}{v_i}$ takes any specific value with probability zero.)
Applying an invertible linear transformation on $V$ preserves properties such as one vector being in
the span of another set. So we can assume w.l.o.g. that the first coordinate in all elements of $V$ is non-zero. Scaling each
$v_i$ by a non-zero scalar also preserves the properties of spans and so we can assume w.l.o.g. that the first coordinate in each
$v_i$ is equal to $1$. Now, for $v_i$ to be in the span of $v_j,v_k$ it must be that either $v_i \in \{v_j,v_k\}$ or $v_i$ is on
the line passing through $v_j,v_k$ (and they are all distinct). Thus, we have a $\delta$-LCC configuration with dimension $n-1$.
\end{proof}

In view of this lemma, in order to prove Theorem~\ref{thm-lcc} it is enough to prove:

\begin{thm}\label{thm-lccconf}
Let $V = (v_1,\ldots,v_m)\in (\C^d)^m$ be a $\delta$-LCC configuration. Then
\[ \dimension(V) \leq O(1/\delta^9). \]
\end{thm}

\subsection{Proof of Theorem~\ref{thm-lccconf}}

Let $V=(v_1,\ldots,v_m)$ be the list of $m$ points in $\C^d$. The main difficulty in proving the theorem is that some of these
points may be the same. That is, two points $v_i,v_j$ can actually correspond to the same vector in $\C^d$. In this case we say that
$v_i,v_j$ are \emph{copies} of each other. Otherwise, we say that $v_i,v_j$ are \emph{distinct}. If $v$ is a point in the list $V$, we let the \emph{multiplicity} of $v$, denoted $M(v)$, be the number of times that (a copy of) $v$ occurs in $V$.

We note that while repetitions make the proof of Theorem~\ref{thm-lccconf} more complicated, we do not know if they actually help
in constructing LCCs with better parameters. Our proof will proceed in an iterative way, at each step identifying a sufficiently
large sublist with small dimension and removing it. The key step will be the following theorem:

\begin{thm}\label{thm-lowranksub}
There exists an integer $K_1 > 0$ s.t. the following holds.  Let $V = (v_1,\ldots,v_m) \in (\C^d)^m$ be a $\delta$-LCC
configuration. Then there exists a sublist $V' \subset V$ of size at least $\delta^3 m/K_1$ and dimension at most $K_1/\delta^6$.
\end{thm}
\begin{proof}
If there exists a point $v \in V$ with multiplicity larger that $\delta m/10$ then the theorem is true
by taking $V'$ to be all copies of this point.
This avoids the case where a point is recovered mostly by its own copies.
For the rest of the proof we can, thus, assume the following.
\begin{fact}\label{fact: v many triples}
For all $v \in V$ and for every sublist $\Delta$ of $V$ of size at most $\delta m /2$ there is a collinear triple containing $v$ such that the other two points in the triple are not in $\Delta$ (and are distinct from $v$).
\end{fact}

We will describe a (probabilistic) construction of a family of collinear triples and build a design matrix from it.
We call a triple of points in $V$ \emph{good} if it contains three distinct collinear points.
We define a family $T$ of good triples as follows: For every line $\ell$ that has at least three distinct points in $V$ we will define (randomly) a family $T_\ell$ of good triples (later we will fix the randomness).
The family $T$ will be the union of all these sets.

\begin{remark}
The construction of $T$ we present is probabilistic.
It is possible to construct $T$ explicitly and achieve similar properties.
We choose to present the probabilistic construction
as it is simpler and less technical.
\end{remark}

Let $\ell$ be such a line with $r$ points on it (counting multiplicities).
Denote by $V(\ell)$ the sublist of $V$ containing all points that lie on $\ell$.
We first take the family $F$ of triples on $[r]$ given by Lemma~\ref{lem-triples} and then pick a random one-to-one mapping $\rho : [r] \mapsto V(\ell)$.
For a triple $t$ in $F$ we denote by $\rho(t)$ the triple of points in $V(\ell)$ that is the image
of $t$ under $\rho$.
We take $T_\ell$ to be the set of all triples $\rho(t)$ with $t \in F$ and such that $\rho(t)$ is good (i.e., it `hits' three distinct points).

Intuitively, we will have many good triples on a line (in expectation) if there are no two points whose copies cover most of the line (then the probability of hitting three distinct points is small). We will later show that this cannot happen on too many lines.

The next proposition shows that there is a way to fix the randomness so that $T$ contains a quadratic number of triples.
\begin{prop}\label{prop-expectation}
The expectation of $|T|$ is at least $\alpha m^2$
with $\alpha = (\delta / 15)^3$.
\end{prop}
We will prove this proposition later in Section~\ref{sec-proofprop} and will continue now with the proof of the theorem.

Fix $T$ to be a family of triples that has size at least the expectation of $|T|$.
By construction and Lemma~\ref{lem-triples}, the family $T$ contains only good triples and each pair of points appears in at most $6$ different triples (since every two distinct points define a single line and two non-distinct points never appear in a triple together).
The family $T$ thus defines a $3$-regular hypergraph with vertex set $[m]$
and at least $\alpha m^2$ edges
and of co-degree at most $6$.
Lemma~\ref{lem: comb lem} thus implies that
there is a sublist $V'$ of $V$ of size at least
$$|V'| = m' \geq \alpha m / 12 \geq (\delta/45)^3 m$$ with the following property:
Let $T'$ be the subfamily of $T$ that $V'$ induces.
Every $v'$ in $V'$ is contained in at least $\alpha m /2$ triples in $T'$.

By a slight abuse of notation,
we also denote by $V'$ the $m' \times d$ matrix with rows defined by the points in $V'$ (including repetitions). We now use the triples in $T'$ to construct a matrix $A'$ so that $A' \cdot V' = 0$. By the above discussion $A'$ is a $(3, \alpha m/2, 6)$-design matrix and thus, by Theorem~\ref{thm-rankdesign}, has rank at least $$ m' - \left( \frac{18 m'}{\alpha m}\right)^2 \geq m' -
(18/\alpha)^2$$ and so $$ \dimension(V') \leq (18/ \alpha)^2 \leq (60/\delta)^6$$
as was required.
\end{proof}

The next proposition shows how the above theorem can be used repeatedly on a given LCC.
\begin{prop}\label{prop-repeat}
There exist an integer $K_2 >0 $ s.t. the following holds:
Let $V = (v_1,\ldots,v_m)\in (\C^d)^m$ be a $\delta$-LCC configuration
and let $U,W$ be a partition of $V$ into two disjoint sublists such that $W \cap \span(U) = \emptyset$. Then there exists a new partition of $V$ to two sublists $U'$ and $W'$ such that $W' \cap \span(U') = \emptyset$ and such that
\begin{enumerate}
\item $|U'| \geq |U| + \delta^3 m / K_2$, and
\item $\dimension(U') \leq \dimension(U) + K_2/\delta^6$.
\end{enumerate}
\end{prop}
\begin{proof}
First, we can assume that all points in $W$ have multiplicity at most $\delta m /2$ (otherwise we can add one point from $W$ with high multiplicity to $U$ to get $U'$). Thus, for all points $v$ and all sublists $\Delta$ of size at most $\delta m /2$ there is a collinear triple of three distinct points containing $v$ and two other points outside $\Delta$. Again, this is to avoid points that are recovered mostly by copies of themselves.

For a point $w \in W$ we define three disjoint sublists of points
$U(w), P_1(w)$ and $P_2(w)$. The first list, $U(w)$, will be
the list of all points in $U$ that are on special lines through $w$ (that is, lines containing $w$ and at least two other distinct points). Notice that, since $w \not\in \span(U)$, each line through $w$ can contain at most one point from $U$. The second list, $P_1(w)$, will be the list of points in $W \setminus \{w\}$ that are on a line containing $w$ and a point from $U$.
The third list, $P_2(w)$, will be of all other points on special lines through $w$ (that is, on special lines that do not
intersect $U$). These three lists are indeed disjoints, since $w$ is the only common point between two lines passing through it.
By the above discussion we have that $|P_1(w)| + |P_2(w)| \geq \delta m /2$ for all $w \in W$ (since removing these two lists
destroys all collinear triples with $w$). We now separate the proof into two cases:

\paragraph{Case I : There exists $w \in W$ with  $|P_1(w)| > \delta m /4$.}
In this case we can simply take $U'$ to be the points in $V$ that are also in the span of $\{w\} \cup U$.
This new $U'$ will include all points in $P_1(w)$ and so
will grow by at least $\delta m /4$ points. Its dimension will grow by at most one and so we are done.

\paragraph{Case II : For all $w \in W$, $|P_2(w)| \geq \delta m /4$.}
Denote $m' = |W|$.
In this case $W$ itself is a $\delta'$-LCC configuration with $$\delta' = \frac{\delta m}{8m'}.$$
Applying Theorem~\ref{thm-lowranksub} we get a sublist $U'' \subset W$ of size at least
$$ \frac{(\delta')^3 m'}{K_1}
\geq (\delta/8)^3 \cdot \frac{m}{K_1}$$ and dimension at most
$$\frac{K_1}{(\delta')^{6}} \leq K_1 (8/\delta)^6.$$
We can thus take $U'$ to be the points in $V$ that are in the span of $U \cup U''$ and the proposition is proved.
\end{proof}

\begin{proof}[Proof of Theorem~\ref{thm-lccconf}]
We apply Proposition~\ref{prop-repeat} on $V$, starting with the partition $U = \emptyset, W = V$ and ending when $U = V, W = \emptyset$. We can apply the proposition at most $K_2/\delta^3$
times and in each step add at most $K_2/\delta^6$ to the dimension of $A$ (which is initially zero).
Therefore, the final list $U = V$ will have dimension at most $O(1/\delta^{9})$.
\end{proof}

\subsection{Proof of Proposition~\ref{prop-expectation}}\label{sec-proofprop}

Order the points in $V$ so that all copies of the same point are consecutive  and so that $M(v_i) \leq M(v_j)$ whenever $i \leq j$. Let $S \subset V$ be the sublist containing the first $\delta m /10$ points in this ordering (we may be splitting the copies of a single point in the middle but this is fine). We will use the following simple fact later on:
\begin{fact}\label{fact-S}
If $v \in S$ and $M(v') < M(v)$ then $v' \in S$.
\end{fact}

For a point $v \in V$ we denote by $T(v)$ the set of (ordered) triples in $T$ containing $v$ and
for a line $\ell$ by $T_\ell(v)$ the set of (ordered) triples in $T_\ell$ containing $v$.
Recall that these are all random variables determined by the choice of the mappings $\rho$ for each line $\ell$.

The proposition will follow by the following lemma.

\begin{lem}
\label{lem: T(v) large}
Let $v \in S$. Then the expectation of $|T(v)|$ is at least $(\delta / 10)^2 m$.
\end{lem}

The lemma completes the proof of the proposition:
summing over all points in $S$ we get
\begin{eqnarray*}
\E[|T|] &\geq& \E\left[ (1/3)\sum_{v \in V} |T(v)| \right]\,\,\,\,\, \text{\small (each triple is counted at most three times)}\\
& \geq& (1/3)\sum_{v \in S} \E[ |T(v)| ] \\
&\geq& (1/3) \cdot (\delta m /10) \cdot ( (\delta / 10)^2 m ) \geq
(\delta / 15)^3 m^2 .
\end{eqnarray*}

\begin{proof}[Proof of Lemma~\ref{lem: T(v) large}]

Denote by $L(v)$ the set of all special lines through $v$.
To prove the lemma we will identify a subfamily $L'(v)$ of $L(v)$
that contributes many triples to $T(v)$.
To do so, we need the following definitions.
For a set $\gamma \subset \C^d$
denote by $P(\gamma)$ the set of distinct points in $V$ that are in $\gamma$.
Denote $M(\gamma) = \sum_{v \in P(\gamma)} M(v)$.
Denote by $P(\bar S)$ the set of distinct points not in $S$.

\begin{define}[Degenerate line]
Let $\ell \in L(v)$.
We say that $\ell \in L(v)$ is {\em degenerate} if either
\begin{enumerate}
\item The size of $P(\ell) \cap P(\bar S)$ is at most one. That is, $\ell$ contains at most one distinct point  outside $S$.
    Or,
\item There exists a point $v_\ell \in P(\ell)$, distinct from $v$,
such that $M(v_\ell) \geq (1 - \delta/10)M(\ell)$.
\end{enumerate}
A degenerate line satisfying the first (second) property above will be called a degenerate line of the first (second) kind.
\end{define}

Define $L'(v)$ as the set of line $\ell$ in $L(v)$ that are not degenerate.
We will continue by proving two claims.
The first claim shows that every line in $L'(v)$ contributes many triples in expectation to  $T(v)$.

\begin{claim}\label{cla-expline}
For every $\ell \in L'(v)$ we have $\E[ |T_\ell(v)| ] \geq \delta M(\ell) /10$.
\end{claim}

\begin{proof}
Denote $r = |M(\ell)|$. The family of triples $T_\ell$ is obtained by taking a family of $r(r-1)$ triples $F$ on $[r]$ (obtained
from Lemma~\ref{lem-triples}) and mapping it randomly to $\ell$, omitting all triples that are not good (those that do not have
three distinct points). For each triple $t \in F$ the probability that $\rho(t)$ will be in $T_\ell(v)$ can be lower bounded by
$$ \frac{3}{r} \cdot \frac{2r}{3(r-1)} \cdot  \frac{\delta}{20}  = \frac{\delta}{10(r-1)}$$

The factor of $3/r$ comes from the probability that one of the three entries in $t$ maps to $v$ (these are disjoint events so we
can sum their probabilities).

The next factor, $2r/(3(r-1))$, comes from the probability that the second entry in $t$ (in some fixed order) maps to a point
distinct from $v$. Indeed since $|P(\ell) \cap P(\bar S)| \geq 2$ and using Fact~\ref{fact-S} we know that there are at least two
distinct points $v',v''$ on $\ell$ with $M(v') \geq M(v)$ and $M(v'') \geq v$. Since $M(v) + M(v') + M(v'') \leq r$, we get that
 $M(v) \leq r/3$, and so there are at least $2r/3$ ``good'' places for the second point to map to.

The last factor, $\delta / 20$, comes from the probability that the third element of the triple will map to a point distinct from
the first two. The bound of $\delta/20$ will follow from the fact that $\ell$ does not satisfy the second property in the
definition of a degenerate line. To see why, let $v_2$ be the image of the second entry in $t$. Since $\ell$ is not degenerate,
$r' \triangleq  r - M(v_2) > \delta r / 10$. Since $|P(\ell) \cap P(\bar S)| \geq 2$, there is a point $v'$ in $P(\bar S)$ not in
$\{v,v_2\}$, and hence, by Fact~\ref{fact-S}, $M(v) \leq M(v')$. Since $M(v) + M(v') \leq r'$, we get that $M(v) \leq r'/2$. Thus
$r' - M(v) \geq r'/2 \geq \delta r/20$. But $r'-M(v)$ is exactly the number of `good' places that the third entry can map to that
are from $v$ and $v_2$.

Using linearity of expectation we can conclude
$$ \E[ |T_\ell(v)| ] \geq r(r-1) \cdot \frac{\delta}{10(r-1)} = \delta r /10.$$
\end{proof}

The second claim shows that there are many points on lines in $L'(v)$.
\begin{claim}\label{cla-manylines}
With the above notations, we have:
$$ \sum_{\ell \in L'(v)} M(\ell) \geq \delta m /10.$$
\end{claim}
\begin{proof}
Assume in contradiction that
$$\sum_{\ell \in L'(v)} M(\ell) < \delta m / 10.$$
Let $\Delta'$ denote the sub-list of $V$ containing all points that lie on lines in $L'(v)$ so that $|\Delta'| \leq \delta m
/10$. We will derive a contradiction by finding a small sublist $\Delta$ of $V$ (containing $\Delta'$ and two other small
sub-lists) that would violate Fact~\ref{fact: v many triples}. That is, if we remove $\Delta$ from $V$, we destroy all collinear
triples containing $v$.

Let $\ell$ be a degenerate line of the second kind. Then there is a point $v_\ell$ on it that is distinct from $v$ and has
multiplicity at least $(1 - \delta/10)M(\ell)$. For every such line let $\Delta_\ell$ denote the sublist of $V$ containing all of
the at most $(\delta/10)M(\ell)-M(v)$ points on this line that are distinct from both $v$ and  $v_\ell$. Let $\Delta_2$ denote
the union of these lists $\Delta_\ell$ over all degenerate lines of the second kind. We now have that $|\Delta_2| \leq \delta m
/10$ since $\sum_{\ell} (M(\ell)-M(v)) \leq m$ and  in each line $\ell$ we have
$$ |\Delta_\ell| \leq (\delta/10)M(\ell)-M(v) \leq (\delta/10)(M(\ell) - M(v)).$$
Notice that, removing the points in $\Delta_2$ destroys all collinear
triples on degenerate lines of the second kind.

Finally, let $\Delta_S$ denote the sublist of $V$ containing all points that have a copy in  $S$. Thus $\Delta_S$ contains the
list $S$ (of at most $\delta m/10$ elements), plus all of the at most $\delta m/10$ copies of the last point in $S$, meaning that
$|\Delta_S| \leq \delta m/5$. Removing $\Delta_S$ destroys all collinear triples on degenerate lines of the first kind. Define $\Delta$ as the union of
the three sublists $\Delta',\Delta_2$ and $\Delta_S$. From the above we have that removing $\Delta$ from $V$ destroys all
collinear triples containing $V$ and that $|\Delta| \leq 4 (\delta/10)m < \delta m /2$. This contradicts Fact~\ref{fact: v many
triples}.
\end{proof}
%
%
%

Combining the two claims we get that for all $v \in S$,
$$ \E[ |T(v)| ]\geq \sum_{\ell \in L'(v)} \E[ |T_\ell(v) |] \geq \sum_{\ell \in L'(v)} \delta M(\ell) / 10 \geq
(\delta /  10) \cdot (\delta m / 10) = (\delta / 10)^2 m .$$
This completes the proof of Lemma~\ref{lem: T(v) large}.
\end{proof}

\section{Extensions to Other Fields}\label{sec-finite}
In this section we show that our results can be extended from the complex field to fields of characteristic
zero, and even to fields with very large positive characteristic. The argument is quite generic and relies on Hilbert's Nullstellensatz.
\begin{define}[$T$-matrix]
Let $m,n$ be integers and let $T \subset [m] \times [n]$. We call an $m \times n$ matrix $A$ a {\em $T$-matrix} if all entries of $A$ with indices in $T$ are non-zero and all entries with indices outside $T$ are zero.
\end{define}

\begin{thm}[Effective Hilbert's Nullstellensatz \cite{Kol88}]
Let $g_1,\ldots,g_s \in \Z[y_1,\ldots,y_t]$ be degree $d$ polynomials with coefficients in $\{0,1\}$  and let
$$Z \triangleq \{ y \in \C^t \,|\, g_i(y)=0 \,\, \forall i \in [s] \}.$$
Suppose $h \in \Z[z_1,\ldots,z_t]$ is another polynomial with coefficients in $\{0,1\}$ which vanishes on $Z$.
Then there exist positive integers $p,q$ and polynomials $f_1, \ldots, f_s \in \Z[y_1,\ldots,y_t]$ such that
\[ \sum_{i=1}^s f_i \cdot g_i \equiv p \cdot h^q. \]
Furthermore, one can bound $p$ and the maximal absolute value of the coefficients of the $f_i$'s by an explicit function $H_0(d,t,s)$.
\end{thm}

\begin{thm}\label{thm-finitefield}
Let $m,n,r$ be integers and let $T \subset [m] \times [n]$. Suppose that all complex $T$-matrices have rank at least $r$. Let $\F$ be a field of either characteristic zero or of finite large enough characteristic $p > P_0(n,m)$, where $P_0$ is some explicit function of $n$ and $m$. Then, the rank of all $T$-matrices over $\F$  is at least $r$.
\end{thm}
\begin{proof}
Let $g_1,\ldots,g_s \in \C[\{x_{ij} \ | \  i \in [m], j \in [n]\} ]$ be the determinants of all $r \times r$ sub-matrices of an $m \times n$ matrix of variables $X = (x_{ij}$). The statement ``all $T$-matrices have rank at least $r$'' can be phrased as ``if $x_{ij}=0$ for all $(i,j) \not\in T$ and $g_k(X)=0$ for all $k \in [s]$ then $\prod_{(i,j) \in T} x_{ij}=0$.'' That is, if all entries outside $T$ are zero and $X$ has rank smaller than $r$ then it must have at least one zero entry also inside $T$. From Nullstellensatz we know that there are integers $\alpha, \lambda > 0$ and polynomials $f_1,\ldots,f_s$ and $h_{ij}, (i,j) \not \in T$, with integer coefficients such that
\begin{equation}\label{eq-rankpoly}
\alpha \cdot \left( \prod_{(i,j) \in T} x_{ij} \right)^\lambda \equiv \sum_{(i,j) \not\in T} x_{ij} \cdot h_{ij}(X)  + \sum_{k=1}^s f_i(X)\cdot g_i(X).
\end{equation}
This identity implies the high rank of $T$-matrices also over any field $\F$ in which $\alpha \neq 0$. Since we have a bound on $\alpha$ in terms of $n$ and $m$ the result follows.
\end{proof}

\section{Discussion and Open Problems}\label{sec-open}

Our rank bound for design matrices has a dependence on $q$, the number of non-zeros in each row. Can this dependency be removed? This  might be possible since a bound on $q$ follows indirectly from specifying the bound on $t$, the sizes of the intersections. Removing this dependency might also enable us to argue about square matrices. Our results so far are interesting only in the range of parameters where the number of rows is much larger than the number of columns.

With respect to Sylvester-Gallai configurations, the most obvious open problem (discussed in the introduction) is to close the gap between our bound of $O(1/\delta^2)$ on the dimension of $\delta$-SG configuration and the trivial lower bound of $\Omega(1/\delta)$ obtained by a simple partition of the points into $1/\delta$ lines.

Another interesting direction is to explore further the connection between design-matrices and LCCs. The most natural way to construct an LCC is by starting with a low-rank design matrix and then defining the code by taking the matrix to be its parity-check matrix. Call such codes {\em design-LCCs}.
Our result on the rank of design matrices shows, essentially, that design-LCCs over the complex numbers cannot have good parameters in general (even for large query complexity). It is natural to ask whether there could exist LCCs that do not originate from designs. Or, more specifically, whether any LCC defines another LCC (with similar parameters) which is a design-LCC. This question was already raised in \cite{BIW07}. Answering this question over the complex numbers will, using our results, give bounds for general LCCs. It is not out of the question to hope for bounds on LCCs with query complexity as large as polynomial in $m$ (the encoding length). This would be enough to derive new results on rigidity via the connection made in \cite{Dvi10}. In particular, our results on design matrices still give meaningful bounds (on design-LCCs) in this range of parameters.

More formally, our results suggest a bound of roughly $\poly(q,1/\delta)$ on the dimension of $(q,\delta)$-LCCs that arise from designs. A strong from of a conjecture from \cite{Dvi10} says that an LCC $C \subset \F^n$ with $q = n^{\eps}$ queries and error $\delta = n^{-\eps}$, for some constant $\eps>0$, cannot have dimension $0.99 \cdot n$. This conjecture, if true, would lead to new results on rigidity. Thus, showing that any LCC defines a design (up to some polynomial loss of parameters), combined with our results, would lead to new results on rigidity.

%

\section*{Acknowledgements}
We thank Moritz Hardt for many helpful conversations. We thank Jozsef Solymosi
for helpful comments.

\bibliographystyle{alpha}

\bibliography{designrank}

\end{document}